\documentclass[12pt]{article}
\usepackage{amssymb}
\usepackage{amsmath}
\usepackage[dvips]{graphicx}
\begin{document}

\title{\LARGE\bf On the Fourier expansion method for highly accurate computation of the
\\ Voigt/complex error function
\\ in a rapid algorithm}

\author{\normalsize\bf S. M. Abrarov\footnote{\scriptsize{Dept. Earth and Space Science and Engineering, York University, Toronto, Canada, M3J 1P3.}}\, and B. M. Quine$^{*}$\footnote{\scriptsize{Dept. Physics and Astronomy, York University, Toronto, Canada, M3J 1P3.}}}

\date{May 7, 2012}
\maketitle

\begin{abstract}
In our recent publication \cite{Abrarov} we presented an exponential series approximation suitable for highly accurate computation of the complex error function in a rapid algorithm. In this Short Communication we describe how a simplified representation of the proposed complex error function approximation makes possible further algorithmic optimization resulting in a considerable computational acceleration without compromise on accuracy.
$$\\$$
\end{abstract}

Consider the exponential approximation of the complex error function, obtained by Fourier expansion method \cite{Abrarov}
\footnotesize
\begin{equation}\label{eq_1}
w\left( z \right) \approx \frac{i}{{2\sqrt \pi  }}\left[ {\sum\limits_{n = 0}^N {{a_n}{\tau _m}\left( {\frac{{1 - {e^{i\,\,\left( {n\pi  + {\tau _m}z} \right)}}}}{{n\,\pi  + {\tau _m}z}} - \frac{{1 - {e^{i\,\,\left( { - n\pi  + {\tau _m}z} \right)}}}}{{n\,\pi  - {\tau _m}z}}} \right)}  - {a_0}\frac{{1 - {e^{i{\tau _m}z}}}}{z}} \right],
\end{equation}
\normalsize
where $z = x + iy$ is a complex argument and 
$$
{a_n} \approx \frac{{2\sqrt \pi  }}{{{\tau _m}}}\exp \left( { - \frac{{{n^2}{\pi ^2}}}{{\tau _m^2}}} \right)
$$
is a set of the Fourier expansion coefficients.

Let us take a term corresponding to index $n = 0$ apart from the summation sign in approximation series \eqref{eq_1}. This leads to the following representation
\footnotesize
\begin{equation}\label{eq_2}
w\left( z \right) \approx i\frac{{1 - {e^{i{\tau _m}z}}}}{{{\tau _m}z}} + \frac{{i{\tau _m}}}{{2\sqrt \pi  }}\sum\limits_{n = 1}^N {{a_n}\left( {\frac{{1 - {e^{in\pi }}{e^{i{\tau _m}z}}}}{{n\,\pi  + {\tau _m}z}} - \frac{{1 - {e^{ - i\,n\pi }}{e^{i{\tau _m}z}}}}{{n\,\pi  - {\tau _m}z}}} \right)}.
\end{equation}
\normalsize
The expression inside parentheses of the approximation \eqref{eq_2} can be rearranged as 
\scriptsize
$$
\frac{{1 - {e^{in\pi }}{e^{i{\tau _m}z}}}}{{n\,\pi  + {\tau _m}z}} - \frac{{1 - {e^{ - i\,n\pi }}{e^{i{\tau _m}z}}}}{{n\,\pi  - {\tau _m}z}} = {e^{i{\tau _m}z}}\frac{{n\,\pi \left( {{e^{in\pi }} - {e^{ - in\pi }}} \right) + {\tau _m}z\left( {{e^{in\pi }} + {e^{ - in\pi }}} \right)}}{{{n^2}\,{\pi ^2} - \tau _m^2{z^2}}} - \frac{{2{\tau _m}z}}{{{n^2}\,{\pi ^2} - \tau _m^2{z^2}}}\,.
$$
\normalsize
Since \footnotesize $\left( {{e^{in\pi }} - {e^{ - in\pi }}} \right)/2i = \sin \left( {n\pi } \right) \equiv 0$ \normalsize and \footnotesize $\left( {{e^{in\pi }} + {e^{ - in\pi }}} \right)/2 = \cos \left( {n\pi } \right) \equiv {\left( { - 1} \right)^n}$\normalsize, we can rewrite it in form
$$
\frac{{1 - {e^{in\pi }}{e^{i{\tau _m}z}}}}{{n\,\pi  + {\tau _m}z}} - \frac{{1 - {e^{ - i\,n\pi }}{e^{i{\tau _m}z}}}}{{n\,\pi  - {\tau _m}z}} = 2{\tau _m}z\frac{{{{\left( { - 1} \right)}^n}{e^{i{\tau _m}z}} - 1}}{{{n^2}\,{\pi ^2} - \tau _m^2{z^2}}}.
$$

Ultimately, applying this relation in approximation \eqref{eq_2} yields
\begin{equation}\label{eq_3}
w\left( z \right) \approx i\frac{{1 - {e^{i{\tau _m}z}}}}{{{\tau _m}z}} + i\frac{{\tau _m^2z}}{{\sqrt \pi  }}\sum\limits_{n = 1}^N {{a_n}\frac{{{{\left( { - 1} \right)}^n}{e^{i{\tau _m}z}} - 1}}{{{n^2}\,{\pi ^2} - \tau _m^2{z^2}}}}.
\end{equation}

Comparing the approximation representations \eqref{eq_1} and \eqref{eq_3} with each other, we can see that the second one is substantially simplified and we can take advantage of it for computational optimization.

In order to describe shortly and conveniently how the simplified representation of the complex error function \eqref{eq_3} can be appropriately implemented in a rapid algorithm, we will use the standard Matlab notations for element-by-element array multiplication  $".*"$ and division $"./"$ with regular and bold letters for scalars and arrays, respectively. Thus the algorithmic implementation based on approximation \eqref{eq_3} can be symbolically conformed to
$$
{\mathbf{w}}: = i\left( {1 - {\mathbf{B}}} \right)./{\mathbf{A}} + i\frac{{{\tau _m}}}{{\sqrt \pi  }}{\mathbf{A}}.*\sum\limits_{n = 1}^N {{a_n}\left[ {{{\left( { - 1} \right)}^n}{\mathbf{B}} - 1} \right]./\left( {{n^2}\,{\pi ^2} - {\mathbf{C}}} \right)},
$$
where ${\mathbf{A}}: = {\tau _m}{\mathbf{z}}$, ${\mathbf{B}}: = {e^{i{\mathbf{A}}}}$, and ${\mathbf{C}}: = {\mathbf{A}}.*{\mathbf{A}}$ are the only three arrays involved in computation (after initialization of ${\mathbf{A}}$, the array ${\mathbf{z}}$ is never used anymore and can be immediately cleaned from the computer memory). Since none of these arrays are dependent on index $n$, each of them are computed just a single time. As follows from this scheme, only the array ${\mathbf{B}}$ requires the exponentiation. All other operations are trivial such as summation, multiplication and division. Since the calculation of the array ${\mathbf{B}}$ is made only once, it has a very minor contribution (at ${\tau _m} = 12$ and $N = 23$ it is less than $6\%$) for the total computation time. This explains why the calculation of the Voigt/complex error function is {\it de facto} as rapid as a rational approximation.

The notation ${e^{i{\mathbf{A}}}}$ implies element-wise exponential operation over the array  ${\mathbf{A}}$ that returns the array ${\mathbf{B}}$ with same dimension; such element-wise functional operations are typical in array programming languages like Matlab/GNU Octave or modern FORTRAN. 

It should be noted that if the high-accuracy is not needed, the acceleration can be additionally gained by reducing the values of ${\tau _m}$ and $N$, say to $9$ and $12$, respectively. The approximation of the exponential function ${e^{i{\tau _m}z}}$ by a rational function is not required as its contribution to the total computation time remains insignificant even in this case.

The advantage of the approximation \eqref{eq_3} implementation becomes particularly evident at an extended size of the input array. Specifically, the computational test reveals that with size of input array ${\mathbf{z}}$ exceeding $30$ million elements, the algorithm based on approximation \eqref{eq_3} is more than twice faster than the one based on approximation \eqref{eq_1} and, compared to the Weideman$'$s algorithm (we used Matlab code in page 1511, Ref. \cite{Weideman}), it becomes faster by a factor greater than $5$ (see also Ref. \cite{Abrarov} for comparison).

Similar to approximation \eqref{eq_1} the numbers generated by approximation \eqref{eq_3} match up to the last decimal digits with those generated by well-known Algorithm 680 \cite{Abrarov, Poppe}. This signifies that due to simplified form of the approximation \eqref{eq_3} the considerable acceleration in computation is achieved without compromise on accuracy.

All calculations were performed with help of Matlab 7.9.0 on a typical desktop computer Intel(R) Quad CPU with RAM 4.00 GB.


------------
\\ {\footnotesize Submitted to journal on October 18, 2011.}

\end{document}